# Classical Theory of Fourier Series : Demystified and Generalised

VIVEK V. RANE

The Institute of Science,

15, Madam Cama Road,

Mumbai-400 032

e-mail address : v_v_rane@yahoo.co.in

**Abstract** :  For a Riemann integrable function on an interval and for a point therein, we define 'Fourier Series at the point on the interval' and bring out how and when the function element becomes expressible as a Fourier series . In this process, we also generalise the Fourier theory ,bringing in such concepts as finite Fourier series, local Fourier series , right/left hand Fourier series. For an integer $q > 1$ , we also sum up  the typical   subseries modulo q , of the basic Fourier series.

# Classical Theory of Fourier series : Demystified and Generalised

VIVEK V. RANE

DEPARTMENT OF MATHEMATICS,

THE INSTITUTE OF SCIENCE,

15, MADAM CAMA ROAD,

MUMBAI-400 032

INDIA

e-mail address : v_v_rane@yahoo.co.in

Given a finite interval , a point in it and a Riemann integrable function on the interval, we formally define 'trigonometric (exponential) Fourier series of the function , at the point on the interval' and we bring out in all its generalities how and when the function element becomes expressible as the Fourier series . We also give an immediate proof of the Cesa'ro or (C,1) summability of the Fourier series. We shall need a good deal of Riemann-Stieltjes integration theory, for which we refer the reader to the book of Tom Apostol [1]. All this will be done using the two facts namely

1) generalized Euler's summation formula

2) the boundedly convergent series expression $u - [u] - \frac{1}{2} = -\sum_{|n| \geq 1} \frac{e^{2\pi i n u}}{2\pi i n}$

for non-integral u , where [u] is the integral part of the real variable u .

This approach enables us in author [ 4 ] to identify and study the Fourier series of the derivatives (with respect to first variable) of Hurwitz zeta function $\zeta(s,\alpha)$, as a function of second variable $\alpha$ on the unit interval [0,1]. There, we deal with finite Fourier series of $\zeta(s,\alpha)$ for a rational value of $\alpha$ in the unit interval. There we also deal with Lerch's zeta function on similar lines. This approach was first resorted to, in author[3]. In the present context, we also refer the reader to author [ 2 ].

Let $f$ be a function defined and Riemann integrable on an interval $I = [c,d]$ of length $2\delta = d - c$ and let x be a point in $[c,d]$. We define the Fourier series of $f$ at the point $x$ on the interval $[c,d]$ and write

$$f(x) \sim \sum_{n=-\infty}^{\infty} c_n e^{\frac{\pi i n x}{\delta}}, \text{ if } c_n = \frac{1}{2\delta}\int_e^d f(u) e^{-\frac{\pi i n u}{\delta}}. \text{ (Exponential Fourier Series)}$$

and we write $f(x) \sim \frac{a_0}{2} + \sum_{n \geq 1}\left(a_n \cos \frac{\pi n x}{\delta} + b_n \sin \frac{\pi n x}{\delta}\right)$,

if $a_n = \frac{1}{2\delta}\int_c^d f(u)\cos \frac{\pi n u}{\delta} du$ and $b_n = \frac{1}{2\delta}\int_c^d f(u)\sin \frac{\pi n u}{\delta} du$. (Trigonometric Fourier series).

Notation : For a real variable u, we write $[u]$ for its integral part. We write $\psi(u) = u - [u] - \frac{1}{2}$. Note that $\psi(u) = -\sum_{|n| \geq 1} \frac{e^{2\pi i n u}}{2\pi i n}$, if u is not an integer. For x with $c \leq x < d$, we write $f(x+) = \lim_{\epsilon \to 0} f(x+\epsilon)$, where $\epsilon > 0$, if this limit exists.

Similarly, for x with $c < x \leq d$, we write $f(x-) = \lim_{\epsilon \to 0} f(x-\epsilon)$ where $\epsilon > 0$, if the limit exists. We call $f(x+)$, the right hand limit of $f$ at $x$. Similarly, we define $f(x-)$, the left hand limit of $f$ at the point $x$.

If the sequence of arithmetic means of the partial sums of a given series converges, then the given series is said to be *Cesàro* summable or $(C,1)$ summable. Using our approach, we can generalize the concept of Fourier series to local Fourier series, local right hand/left hand Fourier series and finite Fourier series. In fact, we prove the following.

**Theorem 1** : Let f be a function defined and Riemann integrable on the interval $[c,d]$ of length $2\delta = d - c$ and let $\psi(u) = u - [u] - \frac{1}{2}$. Let a be a point in $(c,d)$.

I) If $f(a-)$ exists, then we have

$$\frac{f(a-)}{2} = \frac{1}{2\delta}\int_c^a f(u)du + \int_c^a \psi\left(\frac{u-a}{2\delta}\right)df(u) + f(c)\left(\frac{1}{2} - \frac{a-c}{2\delta}\right),$$

provided $f(a) = f(a-)$.

II) If $f(a+)$ exists, then we have

$$\frac{f(a+)}{2} = \frac{1}{2\delta}\int_a^d f(u)du + \int_a^d \psi\left(\frac{u-a}{2\delta}\right)df(u) + f(d)\left(\frac{1}{2} - \frac{d-a}{2\delta}\right),$$

provided $f(a) = f(a+)$.

III) If both $f(a+)$ and $f(a-)$ exist and if $f(c) = f(d)$, then we have

$$\frac{f(a+)+f(a-)}{2} = \frac{1}{2\delta}\int_c^d f(u)du + \int_c^d \psi\left(\frac{u-a}{2\delta}\right)df(u).$$

In addition, if we can interchange $\int_c^d$ and $\sum_n$ in the Riemann-Stieltjes integral

$\int_c^d \psi\left(\frac{u-a}{2\delta}\right) df(u)$, that is, if $\int_c^d \left( \sum_{|n|\geq 1} \frac{e^{2\pi i n\left(\frac{u-a}{2\delta}\right)}}{2\pi i n} \right) df = \sum_{|n|\geq 1} \frac{1}{2\pi i n} \int_c^d e^{2\pi i n\left(\frac{u-a}{2\delta}\right)} df$, then we have

the local Fourier series $\frac{f(a+)+f(a-)}{2} = \sum_{n=-\infty}^{\infty} c_n e^{\frac{\pi i n a}{\delta}}$, where $c_n = \frac{1}{2\delta} \int_c^d f(u) e^{-\frac{\pi i n u}{\delta}} du$.

Otherwise, the local Fourier series converges to $\frac{f(a+)+f(a-)}{2}$ in (c,1) sense.

If f is periodic with period d-c=$2\delta$, then $c_n = \int_{-\delta}^{\delta} f(u) e^{-\frac{\pi i n u}{\delta}} du = \int_0^{2\delta} f(u) e^{-\frac{\pi i n u}{\delta}} du$ and we get the usual Fourier series.

**Note** : For real $\alpha, \beta$ if

$\int_c^d \psi(\alpha u + \beta) df(u) = -\int_c^d \sum_{|n|\geq 1} \frac{e^{2\pi i n(\alpha u + \beta)}}{2\pi i n} df(u) = -\sum_{|n|\geq 1} \frac{1}{2\pi i n} \int_c^d e^{2\pi i n(\alpha u + \beta)} df(u)$,

then the function f will be called a 'good function on the interval $[c,d]$'.

If $f$ is a function of bounded variation on $[c,d]$, then the Riemann – Stieltjes integral behaves like Riemann integral and thus $\sum_n$ and $\int_c^d$ can be interchanged. On the other hand, if $f$ is differentiable with its derivative $f'$ Riemann or Lebesgue integrable on $I = [c,d]$, then the Riemann-Stieltjes integral $\int_c^d \psi\left(\frac{u-a}{2\delta}\right) df$ becomes Riemann integral $\int_c^d \psi\left(\frac{u-a}{2\delta}\right) f'(u) du$ and thus $\sum_n$ and $\int_c^d$ are interchangeable.

**Corollary** : I) Let $f$ be Riemann integrable on an interval $[a-\delta, a]$ and let

$f(a-)$ exist. Then we have $\frac{f(a-)}{2} = \frac{1}{2\delta} \int_{a-\delta}^{a} f(u)du + \int_{a-\delta}^{a} \psi\left(\frac{u-a}{2\delta}\right)df$

and if $f$ is good on the interval $(a-\delta, a)$, then we have **left hand local Fourier**

**series** namely, $f(a-) = \frac{1}{2\delta} \sum_{n=-\infty}^{\infty} e^{\frac{\pi i n a}{\delta}} \int_{a-\delta}^{a} f(u) e^{-\frac{\pi i n u}{\delta}} du$.

II)   Let $f$ be Riemann integrable on the interval $[a, a+\delta]$ and let $f(a+)$ exist.

Then $\frac{f(a+)}{2} = \frac{1}{2\delta} \int_{a}^{a+\delta} f(u)du + \int_{a}^{a+\delta} \psi\left(\frac{u-a}{2\delta}\right)df$.

In addition, if $f$ is good on $[a, a+\delta]$, we have **right hand local Fourier series**

namely, $f(a+) = \frac{1}{2\delta} \sum_{n=-\infty}^{\infty} e^{\frac{\pi i n a}{\delta}} \int_{a}^{a+\delta} f(u) e^{-\frac{\pi i n u}{\delta}} du$.

III) Let $f$ be be Riemann intergrable on the interval $[a-\delta, a+\delta]$ and let both

$f(a-)$ and $f(a+)$ exist. Then, we have $\frac{f(a+)+f(a-)}{2} = \frac{1}{2\delta} \int_{a-\delta}^{a+\delta} f(u)du + \int_{a-\delta}^{a+\delta} \psi\left(\frac{u-a}{2\delta}\right)df$.

In addition, if $f$ is good on $[a-\delta, a+\delta]$, we have the **local Fourier series**

namely

$\frac{f(a+)+f(a-)}{2} = \frac{1}{2\delta} \sum_{n=-\infty}^{\infty} e^{\frac{\pi i n a}{\delta}} \int_{a-\delta}^{a+\delta} f(u) e^{-\frac{\pi i n u}{\delta}} du$.

Before stating our next Theorems, we introduce additional notations.

**Additional Notations** : For a given integer $q \geq 1$ and a given integer r with

$1 \leq r < q$, consider the function $\psi(u, r, q) = - \sum_{\substack{|n| \geq 1 \\ n \equiv r (\mod q)}} \frac{e^{2\pi i n u}}{2\pi i n}$. We shall show, in what

follows, that $\psi(u, r, q)$ is a step function for $0 < u < 1$.

If u is an integer, then

$$\psi(u,r,q) = -\sum_{|n|\geq 1} \tfrac{1}{2\pi i n} = -\tfrac{1}{2\pi i} \sum_{m=-\infty}^{\infty} \tfrac{1}{mq+r} = -\tfrac{1}{2\pi i} \{\tfrac{1}{r} + \sum_{m\geq 1}(\tfrac{1}{mq+r} + \tfrac{1}{r-mq})\}$$

$$= -\tfrac{1}{2\pi i}\{\tfrac{1}{r} + \sum_{m\geq 1}\tfrac{2r}{r^2-m^2q^2}\} = -\tfrac{1}{2\pi i}\left(\tfrac{1}{r} - 2r\sum_{m\geq 1}\tfrac{1}{m^2q^2-r^2}\right) = c(r,q), \text{ say.}$$

Note that $\psi(u,0,q) = -\sum_{\substack{|n|\geq 1 \\ n\equiv 0(\bmod q)}} \dfrac{e^{2\pi i n u}}{2\pi i n} = \dfrac{1}{q}(qu - [qu] - \dfrac{1}{2})$ for non-integral u.

Next, we state our next Theorems.

**Theorem 2** : Let f be Riemann integrable on the unit interval $[0,1]$ with $f(0) = f(1)$. For an integer $q \geq 1$, let f be continuous at the rational points $\tfrac{\ell}{q}$ for $\ell = 1,2,.....,q$. Let a be an integer with $1 \leq a \leq q$.

Then $f\left(\tfrac{a}{q}\right) = \sum_{r-1}^{q} b_r e^{\tfrac{2\pi i r a}{q}}$ with $b_r = \int_0^1 \psi(u,q-r,q)df(u)$ for $r = 1,2,.....,q-1$ and

$$b_0 = \int_0^1 f(u)du + \int_0^1 \psi(u,0,q)df .$$

**Remark** : Note that for the given integer $q \geq 1$, the coefficients $(b_r)$'s are independent of the integer a, where a varies over the set $\{1,2,3,.....,q\}$. Thus if $x = \tfrac{a}{q}$ is a rational number in the unit interval $[0,1]$, then there are two Fourier series at x, namely 1) the usual infinite Fourier series $f(x) \sim \sum_{n=-\infty}^{\infty} c_n e^{2\pi i n x}$ with

$$c_n = \int_0^1 f(u)e^{-2\pi i n u}du$$ ; and the finite Fourier series modulo q namely,

$$f\left(\tfrac{a}{q}\right) = \sum_{r=1}^{q} b_r e^{\frac{2\pi i r a}{q}}$$. However if $\tfrac{a}{q} = \tfrac{b}{k}$, where $b, k \geq 1$ are integers, then

$$f\left(\tfrac{b}{k}\right) = \sum_{\ell=1}^{k} b'_\ell e^{\frac{2\pi i \ell b}{k}}$$, where $(b'_\ell)_{\ell=1,2,\ldots,k}$ is independent of the integer b, with b

lying in the set $\{1, 2, \ldots, k\}$. Thus, we have a finite Fourier series modulo k.

It is to be noted that the concept of finite Fourier series at a rational point is implicit in the concept of Fourier series of a periodic arithmetical function, in literature.

We shall prove Theorem 2 using our Theorem 3, which we state below.

**Theorem 3** : For $1 \leq r \leq q$, the function $\psi(u, r, q) = -\sum_{\substack{|n| \geq 1 \\ n \equiv r \pmod q}} \frac{e^{2\pi i n u}}{2\pi i n}$ is a periodic step

function with period 1 and for $0 < u < 1$, we have

$$\psi(u, r, q) = c(r, q) - \tfrac{1}{q}\left(\tfrac{1}{2} + \sum_{1 \leq \ell \leq qu}{}' e^{\frac{2\pi i \ell r}{q}}\right)$$

where dash over $\sum_\ell$ indicates that the term corresponding to $\ell = qu$ is to be halved, if u is of the form $\tfrac{R}{q}$ for some integer $R \geq 1$.

Before we give the proofs of our Theorems, we state a few lemmas.

**Lemma 0** : Let $\alpha$ be a step function defined on an interval $[a, b]$ with jump $\alpha_k$ at

$x_k$ where $x_1, x_2, \ldots x_n$ are the simple discontinuities of the function $\alpha$. Let $f$

be a function defined on $[a,b]$ such that not both $f$ and $\alpha$ have discontinuities from left or from right at each $x_k$. Then $\int_a^b f d\alpha$ exists and we have

$$\int_a^b f d\alpha = \sum_{k=1}^{n} f(x_k)\alpha_k,$$

where $\alpha_k = \alpha(x_k+) - \alpha(x_k-)$ if $a < x_k < b$;

if $x_k = a, \alpha_k = \alpha(x_k+) - \alpha(x_k)$;

and if $x_k = b$, then $\alpha_k = \alpha(x_k) - \alpha(x_k-)$.

**Remark** : Lemma 0 is the Theorem 7.11 of the Chapter 7 of Apostol's book [1].

Next, we state our Lemma 1.

**Lemma 1** : (Generalised Euler's summation formula) : Let $f$ be Riemann-integrable on the interval $[a,b]$ such that $f$ is continuous from left at every integer $n$ with $\alpha < n \leq b$. Then

$$\sum_{a < n \leq b} f(n) = \int_a^b f(u) du + \int_a^b (u - [u] - \tfrac{1}{2}) df(u) + f(a)(a - [a] - \tfrac{1}{2}) - f(b)(b - [b] - \tfrac{1}{2}).$$

**Remark** : The proof follows from Lemma 0.

**Lemma 2** : We have, for $\epsilon > 0$ and $\delta > 0$,

1) If $f(a-)$ exists and if $f(a) = f(a-)$, then $f(a-) = \int_{a-\epsilon}^{a} f(u) d\left[\tfrac{u-a}{\delta}\right]$

for $0 < \epsilon < \delta$.

2) If $f(a+)$ exists and if $f(a) = f(a+)$, then $f(a+) = -\int_{a}^{a+\epsilon} f(u) d\left[\tfrac{a-u}{\delta}\right]$

for $0 < \epsilon < \delta$.

3) If $f$ is continuous at $a$, we have $f(a) = \int_{a-\epsilon_1}^{a+\epsilon_2} f(u) d\left[\frac{u-a}{\delta}\right]$ for $0 < \epsilon_1, \epsilon_2 < \delta$.

**Note** : The lemma follows from Lemma 0.

**Lemma 3 (Integration by parts)** : Let the Riemann-Stieltjes integral $\int f d\alpha$ exist on the interval $[a,b]$. Then $\int_a^b \alpha df$ also exists and we have

$$\int_a^b f d\alpha = f(b)\alpha(b) - f(a)\alpha(a) - \int_a^b \alpha df$$

**Proof of Theorem 1** : Let a be the point in the open interval $(c,d)$ and let $f(a-)$ exist. Take $f(a) = f(a-)$. Let $d - c = 2\delta$. Then by Lemma 2, we have

$$f(a-) = \int_c^a f(u) d\left[\frac{u-a}{2\delta}\right].$$

As f is Riemann-integrable on [c,a], we have

$$f(a-) = \int_c^a f(u) d\left[\frac{u-a}{2\delta}\right] = \frac{1}{2\delta} \int_c^a f(u) du + \int_c^a f(u) d\left(\left[\frac{u-a}{2\delta}\right] - \frac{u-a}{2\delta}\right)$$

$$= \frac{1}{2\delta} \int_c^a f(u) du - \int_c^a f(u) d\left(\frac{u-a}{2\delta} - \left[\frac{u-a}{2\delta}\right]\right)$$

$$= \frac{1}{2\delta} \int_c^a f(u) du + \left[f(u)\left(\frac{u-a}{2\delta} - \left[\frac{u-a}{2\delta}\right]\right)\right]_{u=a}^{u=c} + \int_c^a \left(\frac{u-a}{2\delta} - \left[\frac{u-a}{2\delta}\right]\right) df(u), \text{ on integration by parts.}$$

$$= \frac{1}{2\delta}\int_c^a f(u)du + f(c)(\frac{c-a}{2\delta} - [\frac{c-a}{2\delta}]) + \int_c^a (\frac{u-a}{2\delta} - [\frac{u-a}{2\delta}])df(u)$$

$$= \frac{1}{2\delta}\int_c^a f(u)du + \int_c^a \left(\frac{u-a}{2\delta} - [\frac{u-a}{2\delta}] - \frac{1}{2}\right)df + \frac{1}{2}\int_c^a df + f(c)(1 - \frac{a-c}{2\delta}).$$

Writing $\psi(u) = u - [u] - \frac{1}{2}$, we have

$$f(a-) = \frac{1}{2\delta}\int_c^a f(u)du + \int_c^a \psi(\frac{u-a}{2\delta})df + \frac{1}{2}(f(a) - f(c)) + f(c)(1 - \frac{a-c}{2\delta})$$

Thus $\frac{f(a-)}{2} = \frac{1}{2\delta}\int_c^a f(u)du + \int_c^a \psi(\frac{u-a}{2\delta})df + f(c)(\frac{1}{2} - \frac{a-c}{2\delta}).$

Similarly, if $f(a+)$ exists, take $f(a) = f(a+)$.

As $f$ is Riemann-integrable on $[a,d]$,

we have $f(a+) = -\int_a^d f(u)d[\frac{a-u}{2\delta}] = \frac{1}{2\delta}\int_a^d f(u)du + \int_a^d f(u)d(\frac{a-u}{2\delta} - [\frac{a-u}{2\delta}])$

$$= \frac{1}{2\delta}\int_a^d f(u)du + \left\{[f(u)(\frac{a-u}{2\delta} - [\frac{a-u}{2\delta}])]_{u=a}^{u=d} - \int_a^d (\frac{a-u}{2\delta} - [\frac{a-u}{2\delta}])df\right\}, \text{ on integration by parts.}$$

Thus $f(a+) = \frac{1}{2\delta}\int_a^d f(u)du + f(d)(\frac{a-d}{2\delta} - [\frac{a-d}{2\delta}]) - \int_a^d (\frac{a-u}{2\delta} - [\frac{a-u}{2\delta}])df$

$$= \frac{1}{2\delta}\int_a^d f(u)du - \int_a^d (\frac{a-u}{2\delta} - [\frac{a-u}{2\delta}] - \frac{1}{2})df - \frac{1}{2}\int_a^d df + f(d)(1 + \frac{a-d}{2\delta})$$

This gives $\frac{f(a+)}{2} = \frac{1}{2\delta}\int_a^d f(u)du - \int_a^d \psi(\frac{a-u}{2\delta})df + f(d)(\frac{1}{2} + \frac{a-d}{2\delta})$

Thus $\frac{f(a+)}{2} = \frac{1}{2\delta}\int_a^d f(u)du + \int_a^d \psi(\frac{u-a}{2\delta})df + f(d)(\frac{1}{2} + \frac{a-d}{2\delta})$

Adding ,we get

$$\tfrac{f(a+)+f(a-)}{2} = \tfrac{1}{2\delta}\int_c^d f(u)du + \int_c^d \psi\left(\tfrac{u-a}{2\delta}\right)df + f(d)\left(\tfrac{1}{2}+\tfrac{a-d}{2\delta}\right)+ f(c)\left(\tfrac{1}{2}-\tfrac{a-c}{2\delta}\right).$$

If $f(c)=f(d)$, then $f(d)\left(\tfrac{1}{2}+\tfrac{a-d}{2\delta}\right)+ f(c)\left(\tfrac{1}{2}-\tfrac{a-c}{2\delta}\right)= f(c)\left(1+\tfrac{c-d}{2\delta}\right)= f(c)(1-1)=0.$

Thus if $f(c)=f(d)$, we have

$$\tfrac{f(a+)+f(a-)}{2} = \tfrac{1}{2\delta}\int_c^d f(u)df + \int_c^d \psi\left(\tfrac{u-a}{2\delta}\right)df.$$

Note that if f is continuous at a , we have $f(a) = \int_c^d f(u)d[\dfrac{u-a}{2\delta}]$ and vice-versa .

From $f(a) = \int_c^d f(u)d[\dfrac{u-a}{2\delta}]$, we have

$$f(a) = \int_c^d f(u)d([\dfrac{u-a}{2\delta}]-\dfrac{u-a}{2\delta}) + \dfrac{1}{2\delta}\int_c^d f(u)du,$$

as f is Riemann-integrable on [c,d] .

On integration by parts and on noting $f(c)=f(d)$ , we get

$$f(a) = \dfrac{1}{2\delta}\int_c^d f(u)du + \int_c^d \psi(\dfrac{u-a}{2\delta})df(u).$$

This shows that redefining $f(a)$ in terms of $f(a-), f(a+)$ as

$$f(a) = \dfrac{f(a-)+f(a+)}{2}$$ makes f continuous at a .

Note that the value of $\psi(\dfrac{u-a}{2\delta})$ at $u=a$ is immaterial in the integral

$$\int_c^d \psi(\dfrac{u-a}{2\delta})df(u),$$ as f is continuous at $u=a$.

Next, $\int_c^d \psi\left(\frac{u-a}{2\delta}\right)df = \int_c^d\left(-\sum_{n\geq 1}\frac{\sin \pi n\left(\frac{u-a}{\delta}\right)}{\pi n}\right)df = -\sum_{n\geq 1}\frac{1}{\pi n}\int_c^d \sin\frac{\pi n(u-a)}{\delta}df$ ,

if $\sum_{n\geq 1}$ and $\int_c^d$ are interchangeable .

This gives, $\int_c^d \psi\left(\frac{u-a}{2\delta}\right)df = -\sum_{n\geq 1}\left[\frac{\sin \pi n\left(\frac{u-a}{\delta}\right)}{\pi n}\cdot f(u)\right]_{u=c}^{u=d} + \sum_{n\geq 1}\frac{1}{\delta}\int_c^d f(u)\cos \pi n\left(\frac{u-a}{\delta}\right)du$

$= -f(d)\sum_{n\geq 1}\frac{1}{\pi n}\left(\sin \pi n(\frac{d-a}{\delta}) - \sin \pi n(\frac{c-a}{\delta})\right) + \sum_{n\geq 1}\frac{1}{\delta}\int_c^d f(u)\cos \pi n(\frac{u-a}{\delta})du$ .

Next $\sin \pi n(\frac{d-a}{\delta}) - \sin \pi n(\frac{c-a}{\delta}) = 2\sin \pi n(\frac{d-c}{2\delta})\cdot \cos \pi n(\frac{d+c-2a}{2\delta}) = 0$ , as

$d-c = 2\delta$ . Thus $\int_c^d \psi(\frac{u-a}{2\delta})df(u) = \sum_{n\geq 1}\frac{1}{\delta}\int_c^d f(u)\cos \pi n(\frac{u-a}{\delta})du$

$= \frac{1}{2\delta}\sum_{|n|\geq 1}\int_c^d f(u)e^{-\frac{\pi i n(u-a)}{\delta}} = \frac{1}{2\delta}\sum_{|n|\geq 1}e^{\frac{\pi i n a}{\delta}}\int_c^d f(u)e^{-\frac{\pi i n u}{\delta}}du$ .

Thus $\frac{f(a+)+f(a-)}{2} = \frac{1}{2\delta}\sum_{n=-\infty}^{\infty}e^{\frac{\pi i n a}{\delta}}\int_c^d f(u)e^{-\frac{\pi i n u}{\delta}}du$ , provided $f(c) = f(d)$ and we can

interchange $\int_c^d$ and $\sum_n$ in the Riemann-Stieltjes integral $\int_c^d\left(-\sum_{n\geq 1}\frac{1}{\pi n}\sin\frac{\pi n(u-a)}{\delta}\right)df$.

Next, we prove that the sequence of arithmetic means of partial sums of the Fourier series of $f$ at the point a converges to $\frac{f(a+)+f(a-)}{2}$, if both $f(a+)$ and $f(a-)$ exist . However , the proof will not be very rigourous .

Redefine $f(a) = \frac{f(a-)+f(a+)}{2}$ so that f is continuous at $u = a$ and let $f(c) = f(d)$.

Let $S(f,\ell) = \frac{1}{2\delta} \sum_{|n|=1}^{\ell} e^{\frac{\pi i n u}{\delta}} \int_c^d f(u) e^{-\frac{\pi i n u}{\delta}} du$

$= \frac{1}{2\delta} \sum_{|n|=1}^{\ell} \int_c^d f(u) e^{\frac{\pi i n}{\delta}(u-a)} = \frac{1}{\delta} \sum_{n=1}^{\ell} \int_c^d f(u) \cos \frac{\pi n}{\delta}(u-a) du = \sum_{n=1}^{\ell} \int_c^d f(u) d\left(\frac{\sin \frac{\pi n}{\delta}(u-a)}{\pi n}\right)$

$= \sum_{n=1}^{\ell} \left\{ \left[ f(u) \sin \frac{\pi n(u-a)}{\pi n} \right]_{u=c}^{u=d} - \frac{1}{\pi n} \int_c^d \sin \frac{\pi n}{\delta}(u-a) df \right\}$

$= \sum_{n=1}^{\ell} \left\{ \frac{f(d)}{\pi n} \left(\sin \frac{\pi n}{\delta}(d-a) - \sin \frac{\pi n}{\delta}(c-a)\right) \right\} - \sum_{n=1}^{\ell} \frac{1}{\pi n} \int_c^d \sin \frac{\pi n}{\delta}(u-a) df$, as $f(c)=f(d)$

Thus $S(f,\ell) = f(d) \sum_{n=1}^{\ell} \frac{1}{\pi n} \{ 2 \cos \frac{\pi n}{\delta}\left(\frac{d+c}{2} - a\right) \sin \frac{\pi n}{2\delta}(d-c) \} - \sum_{n=1}^{\ell} \frac{1}{\pi n} \int_c^d \sin \frac{\pi n}{\delta}(u-a) df$

Note that $d - c = 2\delta$. This gives $S(f,\ell) = -\sum_{n=1}^{\ell} \frac{1}{\pi n} \int_c^d \sin \frac{\pi n}{\delta}(u-a) df$

$= \int_c^d \left( \sum_{n \le \ell} \frac{\sin \frac{\pi n}{\delta}(u-a)}{-\pi n} \right) df = \int_c^d \psi\left(\frac{u-a}{2\delta}, \ell\right) df$, say, where $\psi(u,\ell) = \sum_{n \le \ell} \frac{\sin 2\pi n u}{-\pi n}$.

Let $\sigma(f,N) = \sum_{1 \le \ell \le N} S(f,\ell) = \int_0^N S(f,\ell) d\ell = \int_0^N d\ell \int_c^d df(u) \psi\left(\frac{u-a}{2\delta}, \ell\right)$

$= \int_c^d df(u) \int_0^N \psi\left(\frac{u-a}{2\delta}, \ell\right) d\ell = \int_c^d df(u) \left( \sum_{\ell=1}^{N} \psi\left(\frac{u-a}{2\delta}, \ell\right) \right)$.

Thus $\lim_{N \to \infty} \frac{\sigma(f,N)}{N} = \int_c^d df(u) \left( \lim_{N \to \infty} \frac{1}{N} \sum_{\ell=1}^{N} \psi\left(\frac{u-a}{2\delta}, \ell\right) \right)$ provided the limit exists.

Since $\lim_{N \to \infty} \frac{1}{N} \sum_{\ell=1}^{N} \psi(u,\ell) = \psi(u)$, we have $\lim_{N \to \infty} \frac{\sigma(f,N)}{N} = \int_c^d \psi\left(\frac{u-a}{2\delta}\right) df(u)$.

Now as $f(c) = f(d)$, we have $\frac{f(a+)+f(a-)}{2} = \frac{1}{2\delta} \int_c^d f(u) du + \int_c^d \psi\left(\frac{u-a}{2\delta}\right) df$.

Thus, we have $\lim\limits_{N\to\infty} \frac{\sigma(f,N)}{N} = \frac{f(a+)+f(a-)}{2} - \frac{1}{2\delta}\int\limits_c^d f(u)du$ .

Thus, the Fourier series of $f$ at the point a in $[c,d]$ converges to $\frac{f(a+)+f(a-)}{2}$ in (c,1) sense , though , the proof is not very rigourous .

**<u>Proof of Theorem 2</u>** : Let $x = \frac{a}{q}$ be a rational number such that $0 < \frac{a}{q} \leq 1$, where $1 \leq a \leq q$ are integers. As $f(0) = f(1)$ and as f is continous at $u = \frac{a}{q}$ , we have

$f\left(\frac{a}{q}\right) = \int\limits_0^1 f(u)du + \int\limits_0^1 \psi\left(u - \frac{a}{q}\right)df(u)$, where the value of $\psi(u - \frac{a}{q})$ at $u = \frac{a}{q}$

is immaterial . Thus $f(\frac{a}{q}) = \int\limits_0^1 f(u)du + \int\limits_0^1 \left(-\sum\limits_{|n|\geq 1} \frac{e^{2\pi in\left(u-\frac{a}{q}\right)}}{2\pi in}\right)df$

$= \int\limits_0^1 f(u)du + \int\limits_0^1 \sum\limits_{r=0}^{q-1} \left(\sum\limits_{\substack{|n|\geq 1 \\ n\equiv r(\bmod q)}} \frac{e^{2\pi in\left(u-\frac{a}{q}\right)}}{2\pi in}\right)df = \int\limits_0^1 f(u)du + \int\limits_0^1 \sum\limits_{r=0}^{q-1} e^{-\frac{2\pi ira}{q}} \left(-\sum\limits_{\substack{|n|\geq 1 \\ n\equiv r(\bmod q)}} \frac{e^{2\pi inu}}{2\pi in}\right)df$

$= \int\limits_0^1 f(u)du + \int\limits_0^1 \sum\limits_{r=0}^{q-1} e^{-\frac{2\pi ira}{q}} \psi(u,r,q)df(u) = \int\limits_0^1 f(u)du + \sum\limits_{r=0}^{q-1} e^{-\frac{2\pi ira}{q}} \int\limits_0^1 \psi(u,r,q)df$

provided $\int\limits_0^1 \psi(u,r,q)df$ exists for each $r = 1,2,.......,q$. However, the statement of

Theorem 3 shows that $\psi(u,r,q)$ is a step function for each $r = 1,2,.......,q$ with simple discontinuities from left at the points $\frac{\ell}{q}$ for $\ell = 1,2,.......,q$, where $f$ is continuous from left and hence $\int\limits_0^1 \psi(u,r,q)df$ exists for each $r = 1,2,.......,q$.

Thus $f\left(\frac{a}{q}\right) = \left(\int_0^1 f(u)du + \int_0^1 \psi(u,o,q)df\right) + \sum_{r=1}^{q-1} e^{-\frac{2\pi i r a}{q}} \int_0^1 \psi(u,r,q)df = b'_0 + \sum_{r=1}^{q-1} b'_r e^{-\frac{2\pi i r a}{q}}$,

say, where $b'_0 = \int_0^1 f(u)du + \int_0^1 \psi(u,o,q)df$. Writing $b'_r = b_{q-r}$ for $r = 1, 2, \ldots, q-1$

and $b'_0 = b_0$, we have $f\left(\frac{a}{q}\right) = \sum_{r=1}^{q} b_r e^{\frac{2\pi i r a}{q}}$.

**Proof of Theorem 3** : Let $a, q$ be integers with $1 \leq a < q$.

Let $c(a,q) = -\frac{1}{2\pi i} \sum_{m=-\infty}^{\infty} -\frac{1}{2\pi i}\left\{\frac{1}{a} + \sum_{m\geq 1}\left(\frac{1}{mq+a} + \frac{1}{a-mq}\right)\right\} = \frac{1}{2\pi i}\left(-\frac{1}{a} + \sum_{m\geq 1} \frac{2a}{m^2 q^2 - a^2}\right)$

Next, let $0 < x < 1$. For a large positive integer M, let $S_{2M+1}(x) = -\sum_{m=-M}^{M} \frac{e^{2\pi i(mq+a)x}}{2\pi i(mq+a)}$

be the sum of the first $(2M+1)$ terms.

Note $\frac{e^{2\pi i(mq+a)x}}{2\pi i(mq+a)} = \int_0^x e^{2\pi i(mq+a)u}\, du + \frac{1}{2\pi i(mq+a)}$.

Thus, $-S_{2M+1}(x) - \sum_{m=-M}^{M} \frac{1}{2\pi i(mq+a)} = \int_0^x du \left(\sum_{m=-M}^{M} e^{2\pi i(mq+a)u}\right)$

$= \int_0^x du \, e^{2\pi i u(a-Mq)}\left(\frac{1-e^{2\pi i qu(2M+1)}}{1-e^{2\pi i qu}}\right) = \int_0^x e^{2\pi i a u}\left(\frac{\sin\left(M+\frac{1}{2}\right)2\pi qu}{\sin \pi qu}\right) du$.

Let $\delta > 0$ be sufficiently small. Let $R$ be the largest integer such that $\frac{R}{q} \leq x$.

If $\frac{R}{q} < x$, then let $I_1 = (0, \delta) \bigcup_{r=1}^{R} \left(\frac{r}{q} - \delta, \frac{r}{q} + \delta\right)$.

If $\frac{R}{q} = x$, then let $I_1 = (0, \delta) \bigcup_{r=1}^{R-1} \left(\frac{r}{q} - \delta, \frac{r}{q} + \delta\right) \bigcup (x - \delta, x)$ so that $I_1$ is disjoint union of intervals.

: 16 :Let $I_2 = I - I_1$. Thus $I_2$ is also disjoint union of intervals.

Thus $\int_0^x e^{2\pi i a u}\left(\frac{\sin(M+\frac{1}{2})2\pi q u}{\sin \pi q u}\right)du = \left(\int_{I_1} + \int_{I_2}\right)e^{2\pi i a u}\left(\frac{\sin(M+\frac{1}{2})2\pi q u}{\sin \pi q u}\right)du$

Note that on $I_2$, the denominator $\sin \pi qu$ is bounded away from zero. Consider

$\int_L$ , where L is one of the non-overlapping sub-intervals, whose union is $I_2$. Let

$\ell_1 < \ell_2$ be the endpoints of L.

Then $\int_L \frac{e^{2\pi i a u}}{\sin \pi qu}\sin(M+\frac{1}{2})2\pi qu\,du = -\int_{\ell_1}^{\ell_2}\frac{e^{\pi i a u}}{\sin \pi qu}d\left(\frac{\cos(M+\frac{1}{2})2\pi qu}{(M+\frac{1}{2})2\pi q}\right)$

$= -\left[\frac{e^{2\pi i a u}}{\sin \pi qu}\cdot\frac{\cos(M+\frac{1}{2})2\pi qu}{(M+\frac{1}{2})2\pi q}\right]_{u=\ell_1}^{\ell_2} + \int_{\ell_1}^{\ell_2}\frac{\cos(M+\frac{1}{2})2\pi qu}{(M+\frac{1}{2})2\pi q}\cdot\frac{d}{du}\left(\frac{e^{2\pi i a u}}{\sin \pi qu}\right)du$

In view of the factor $(M+\frac{1}{2})$ in denominators, both the integrated part and the integral each $\to 0$, as $M \to \infty$.

Next, we evaluate $\int_{I_1}e^{2\pi i a u}\left(\frac{\sin(M+\frac{1}{2})2\pi qu}{\sin \pi qu}\right)du$ as $M \to \infty$.

Hence, consider $\int_{\frac{r}{q}-\delta}^{\frac{r}{q}+\delta} e^{2\pi i a u}\left(\frac{\sin(M+\frac{1}{2})2\pi qu}{\sin \pi qu}\right)du$.

Writing $\frac{r}{q} + v = u$, we have the above integral

$= \int_{-\delta}^{\delta} e^{2\pi i a\left(\frac{r}{q}+v\right)}\left(\frac{\sin(M+\frac{1}{2})2\pi q\left(\frac{r}{q}+v\right)}{\sin \pi q\left(\frac{r}{q}+v\right)}\right)dv = \int_{-\delta}^{\delta} e^{2\pi i a\left(\frac{r}{q}+v\right)}\left(\frac{\sin(M+\frac{1}{2})2\pi q\left(\frac{r}{q}+v\right)}{\sin \pi q\left(\frac{r}{q}+v\right)}\right)dv$

$= e^{\frac{2\pi i a r}{q}}\int_{-\delta}^{\delta} e^{2\pi i a v}\cdot\frac{\sin(M+\frac{1}{2})(2\pi r + 2\pi vq)}{\sin(\pi r + \pi qv)}dv$

$= e^{\frac{2\pi i a r}{q}}\int_{-\delta}^{\delta} e^{2\pi i a u}\cdot\frac{\sin(M+\frac{1}{2})2\pi qu}{\sin \pi qu}du$

$$= e^{\frac{2\pi i a r}{q}} \int_0^\delta \left(e^{2\pi i a u} + e^{-2\pi i a u}\right) \frac{\sin(M+\frac{1}{2})2\pi q u}{\sin \pi q u} du = 2e^{\frac{2\pi i a r}{q}} \int_0^\delta \frac{\cos 2\pi a u}{\sin \pi q u} \cdot \sin(M+\tfrac{1}{2})2\pi q u\, du$$

$$= 2e^{\frac{2\pi i a r}{q}} \int_0^\delta \frac{u \cos 2\pi a u}{\sin \pi q u} \left(\frac{\sin(M+\frac{1}{2})2\pi q u}{u}\right) du .$$

Using the result namely, if a function $g$ is of bounded variation on interval $[0,\delta]$, then $\lim\limits_{\alpha \to \infty} \frac{2}{\pi} \int_0^\delta g(t) \frac{\sin \alpha t}{t} dt = g(0+)$ and in the light of the fact that

$\lim\limits_{u \to 0} \frac{u}{\sin \pi q u} \cos 2\pi a u = \frac{1}{\pi q}$, we have $\lim\limits_{M \to \infty} \int_{\frac{r}{q}-\delta}^{\frac{r}{q}+\delta} e^{2\pi i a u} \left(\frac{\sin(M+\frac{1}{2})2\pi q u}{\sin \pi q u}\right) du = \frac{1}{q} e^{\frac{2\pi i a r}{q}}$

Similarly, we can show $\lim\limits_{M \to \infty} \int_0^\delta e^{2\pi i a u} \left(\frac{\sin(M+\frac{1}{2})2\pi q u}{\sin \pi q u}\right) du = \frac{1}{2q}$

and $\lim\limits_{M \to \infty} \int_{x-\delta}^x e^{2\pi i a u} \left(\frac{\sin(M+\frac{1}{2})2\pi q u}{\sin \pi q u}\right) du = \frac{1}{2q} e^{\frac{2\pi i a r}{q}}$, if $x$ is a fraction of the form $\frac{r}{q}$. Letting

$M \to \infty$, we get $\psi(x,a,q) = \lim\limits_{M \to \infty} S_{2M+1}(x) = -\sum\limits_{m=-\infty}^{\infty} \frac{1}{2\pi i(mq+a)} + \phi(x,a,q)$, where for

$0 < x < 1$, we have $\phi(x,a,q) = -\frac{1}{q}\left(\frac{1}{2} + \sum\limits_{1 \le \ell \le qx}' e^{\frac{2\pi i \ell a}{q}}\right)$. Here dash in $\sum\limits_\ell'$ indicates the last term is to be halved, if $x = \frac{r}{q}$ for some integer $r$.

This completes the Proof of Theorem 3.